%
%
%
\input amstex.tex
\documentstyle{amsppt}
%
%
\def\section#1{\par\bigpagebreak
    \csname subhead\endcsname #1\endsubhead\par\medpagebreak}
\def\subsection#1{\par\medpagebreak
    \csname subsubhead\endcsname #1 \endsubsubhead}
\def\Theorem#1#2{\csname proclaim\endcsname{Theorem #1} #2 
    \endproclaim}
\def\Corollary#1{\csname proclaim\endcsname{Corollary} #1 \endproclaim}
\def\Proposition#1#2{\csname proclaim\endcsname{Proposition #1} #2 
    \endproclaim}
\def\Lemma#1#2{\csname proclaim\endcsname{Lemma #1} #2 \endproclaim}
\def\Remark#1#2{\remark{Remark {\rm #1}} #2 \endremark}
\def\Proof#1{\demo{Proof} #1\qed\enddemo}
\def\Proofof#1#2{\demo{Proof of #1} #2\qed\enddemo}
%
%
\define\ep{\epsilon}
\define\vep{\varepsilon}
\define\dt{\delta}
\define\Dt{\Delta}
\define\ld{\lambda}

\define\Z{{\Bbb Z}}
\define\Q{{\Bbb Q}}

\define\C{{\Bbb C}}
\define\K{{\Bbb K}}
\define\N{{\Bbb N}}
\define\br#1{{\langle{#1}\rangle}}
\define\half{\frac{1}{2}}
\NoRunningHeads
\topmatter
\title\nofrills
A reproducing kernel for nonsymmetric Macdonald polynomials 
\endtitle
\author
Katsuhisa Mimachi
\footnote"$^{\ast1}$"{Department of Mathematics, Kyushu University;  
Hakozaki 33, Fukuoka 812, Japan\hfill\,}
and Masatoshi NOUMI
\footnote"$^{\ast2}$"{Department of Mathematics, Kobe University; 
Rokko, Kobe 657, Japan\hfill\,}
\endauthor
\affil
\endaffil
\address
\flushpar
{\eightpoint
K.M.: Department of Mathematics,
Kyushu University; 
Hakozaki 33, Fukuoka 812, Japan\newline
M.N.: Department of Mathematics, 
Faculty of Science, Kobe University; 
Rokko, Kobe 657, Japan
}
\endaddress
\email 
mimachi\@math.kyushu-u.ac.jp, \quad
noumi\@math.s.kobe-u.ac.jp
\endemail
\abstract
We present a new formula of Cauchy type for 
the nonsymmetric Macdonald polynomials which 
are joint eigenfunctions of $q$-Dunkl operators. 
This gives the explicit formula for a reproducing kernel on the 
polynomial ring of $n$ variables. 
\endabstract
\endtopmatter
\footnote""{\hskip-12pt {\it 1991 Mathematics Subject Classification\,}: 
Primary 33C50; Secondary 33C55,33C80,33D45}
\document
\section{\S0:  Introduction}
In this paper we propose a new formula of Cauchy type
for the nonsymmetric Macdonald polynomials of type $A_{n-1}$.
This can be regarded as an explicit formula for the reproducing 
kernel of a certain scalar product on the polynomial ring of 
$n$ variables. 
A similar result for nonsymmetric Jack polynomials 
was recently given by Sahi \cite{S}.

\par\medpagebreak

The nonsymmetric Macdonald polynomials $E_\ld(x|q,t)$, introduced by 
Macdonald \cite{Ma1}, are characterized as the joint eigenfunctions 
in the polynomial ring of $n$ variables $x=(x_1,\ldots ,x_n)$, 
for  the commuting family of $q$-Dunkl operators. 
(For the precise definition of $E_\ld(x|q,t)$, see Section 1.) 
We define a meromorphic function $E(x;y|q,t)$ 
in $x=(x_1,\ldots ,x_n)$ and $y=(y_1,\ldots ,y_n)$
by 
$$
E(x;y|q,t)=
\prod_{1\le j<i\le n}\frac{(qtx_iy_j;q)_\infty}{(qx_iy_j;q)_\infty}
\prod_{1\le i\le n}   \frac{(qtx_iy_i;q)_\infty}{(x_iy_i;q)_\infty}
\prod_{1\le i<j\le n}\frac{(tx_iy_j;q)_\infty}{(x_iy_j;q)_\infty}.
\tag{0.1}
$$
The main result of this paper is the following. 
\Theorem{A}{
The function $E(x;y|q,t)$ has the following expansion in terms of 
nonsymmetric Macdonald polynomials:
$$
E(x;y|q,t)=\sum_{\ld\in\N^n} a_\ld(q,t) E_\ld(x|q,t) E_\ld(y|q^{-1},t^{-1}). 
\tag{0.2}
$$
For each composition $\ld\in\N^n$, the coefficient $a_\ld(q,t)$ is given by
$$
a_\ld(q,t)=\prod_{s\in \ld}\frac{1-q^{a(s)+1}t^{l(s)+1}}{1-q^{a(s)+1}t^{l(s)}},
\tag{0.3}
$$
where, for each box $s\in\ld$, $a(s)$ and $l(s)$ are the arm-length and the 
generalized leg-length of $s$ in $\ld$. 
}
\noindent
(See Theorems 2.1 and 2.2.)
\par
Assuming that $q$ and $t$ are complex numbers with $0<|q|,|t|<1$, 
we now consider the meromorphic function $K(x;y|q,t)=E(x;y^{-1}|q,t)$  
on the algebraic torus $(\C^\ast)^n\times(\C^\ast)^n$. 
\Theorem{B}{
For each composition $\ld\in \N^n$, we have 
$$
\left(\frac{1}{2\pi\sqrt{-1}}\right)^n\int_{\Bbb T^n}
 K(x;y|q,t) E_\ld(y|q,t)w(y|q,t)
\frac{dy_1\cdots dy_n}{y_1\cdots y_n}=C_\ld.E_\ld(x|q,t)
\tag{0.4}
$$
for all $x=(x_1,\ldots,x_n)\in \Bbb C^n$ with $|x_j|<1$ $(j=1,\ldots,n)$. 
Here ${\Bbb T}^n=\{y=(y_1,\ldots ,y_n)\in(\C^\ast)^n~;~ 
|y_j|=1\ \ (j=1,\ldots ,n)\}$ 
is the $n$-dimensional torus with the standard orientation, and 
$$
w(y|q,t)=\prod_{1\le i<j\le n}
\frac{(y_i/y_j;q)_\infty(qy_j/y_i;q)_\infty}
{(ty_i/y_j;q)_\infty(qty_j/y_i;q)_\infty}.
\tag{0.5}
$$
The constant $C_\ld$ is given by
$$
C_\ld=C_{\ld^+}=\left(\frac{(qt;q)_\infty}{(q;q)_\infty}\right)^n
\prod_{i=1}^n \frac{(q^{\ld^+_i+1}t^{n-i};q)_\infty}
{(q^{\ld^+_i+1}t^{n-i+1};q)_\infty},
\tag{0.6}
$$
where $\ld^+$ is the partition obtained by reordering the parts of $\ld$. 
}
\noindent
(See Theorem 3.2.)
\par\medpagebreak
After preliminaries on nonsymmetric Macdonald polynomials, we will state our 
main results in Section 2.  
In Section 2, we will prove that $E(x;y|q,t)$ has 
an expansion of the form (0.2), 
and reduce the determination of the coefficients $a_\ld(q,t)$ to the case of 
partitions.  
In this paper we will present two ways of determining the 
coefficients $a_\ld(q,t)$ for partitions $\ld$. 
In Section 3, we determine these coefficients in an analytic way by 
asymptotic analysis of $q$-Selberg type integrals similarly as 
in \cite{Mi2}.  
In this proof we will make use of the evaluation of Cherednik's 
scalar product for nonsymmetric Macdonald polynomials. 
Theorem B will also be formulated in Section 3. 
In Section 4, we give an algebraic proof of (0.3) by using the evaluation 
of the nonsymmetric Macdonald polynomials $E_\ld(x;q,t)$ at 
$x=(t^{-n+1},t^{-n+2},\ldots ,1)$ due to Cherednik \cite{C2}. 
This second proof is an extension of the argument of Sahi \cite{S} 
to the $q$-version. 
%

\section{\S1: Nonsymmetric Macdonald polynomials} 
We first recall the definition of 
{\it nonsymmetric Macdonald polynomials}
of type $A_{n-1}$ in the $GL_n$ version. 
Although we follow the presentation by Macdonald \cite{Ma2} 
in principle, 
we use a slightly different convention which is more convenient 
for our purpose. 

\par\medpagebreak

Let ${\K}[x^{\pm1}]$ be the ring of Laurent polynomials 
in $n$ variables $x=(x_1,\ldots ,x_n)$ with coefficients in the field 
$\K=\Q(q,t^\half)$. 
(Although we use this coefficient field for convenience, 
the use of $t^\half$ is {\it not} essential; 
one could work within $\Q(q,t)$ by modifying the argument appropriately. )
We denote by $\tau=(\tau_1,\ldots ,\tau_n)$ the corresponding $q$-shift 
operators.  
For each $i=1,\ldots,n$, the operator $\tau_i$ acts on ${\K}[x^{\pm1}]$ as a 
$\K$-automorphism such that 
$\tau_{i}(x_j)=x_jq^{\delta_{ij}}$ ($j=1,\ldots ,n$). 
The action of the symmetric group $W=\frak{S}_n$ on 
${\K}[x^{\pm1}]$ will be fixed so that 
each $w\in W$ defines the $\Bbb K$-algebra automorphism such that 
$w.x_j=x_{w(j)}$ for $j=1,\ldots,n$. 
The ring ${\K}[x^{\pm1}]$ is identified with the group-ring $\Bbb K[P]$ of 
the integral weight lattice 
$P=\Z\ep_1\oplus\cdots\oplus\Z\ep_n$. 
As usual, we take the symmetric bilinear form $\br{~,~}$  on $P$ such 
that $\br{\ep_i,\ep_j}=\dt_{ij}$ ($1\le i,j\le n$).
For each $\lambda\in P$, we use the notation of multi-indices 
$$
x^\ld=x_1^{\ld_1}\cdots x_n^{\ld_n}, \quad
\tau^\ld=\tau_1^{\ld_1}\cdots \tau_n^{\ld_n}, 
\tag{1.1}
$$
where $\ld_j=\br{\ld,\ep_j}$ ($j=1,\ldots,n$).
The action of the symmetric group $W=\frak S_n$ of degree $n$ on $P$
will be fixed as $w.\ep_i=\ep_{w(i)}$, or equivalently as 
$(w.\ld)_i=\ld_{w^{-1}(i)}$ $(i=1,\ldots ,n)$ for each $w\in W$. 
Note that the commutation relations among the multiplication operators 
$x^\lambda$, the $q$-shift operators $\tau^\mu$ and the permutations $w\in W$ 
are given as follows:
$$
\tau^\mu x^\ld=x^\ld \tau^\mu q^\br{\mu,\ld},\quad
w x^\ld= x^{w.\ld} w,\quad w \tau^\mu= \tau^{w.\mu} w,
\tag{1.2}
$$
for $\ld,\mu\in P$ and $w\in W$. 
We will use the standard notation of the set of positive roots
$$
\Dt^+=\{\ep_i-\ep_j~;~1\le i<j\le n\},
\tag{1.3}
$$
the simple roots $\alpha_i=\ep_i-\ep_{i+1}$ ($i=1,\ldots ,n-1$) 
and the cone of dominant integral weights
$$
P^+=\{\ld\in P~;~ \br{\alpha_i,\ld}\ge0\ \ (i=1,\ldots ,n-1)\}
=\{\ld\in P~;~ \ld_1\ge\cdots\ge\ld_n\}.
\tag{1.4}
$$
We denote the set of {\it compositions} and 
the set of {\it partitions} with length $\le n$ by 
$L=\N\ep_1\oplus\cdots\oplus\N\ep_n\subset P$ and by $L^+=P^+\cap L$,
respectively, where $\N=\{0,1,2,\cdots\}$.

\par\medpagebreak

In what follows, we will make use of the {\it Demazure-Lusztig operators} 
$T_1,\ldots ,T_{n-1}$ defined by 
$$
T_i=t^{\half} +t^{-\half}\frac{1-tx_i/x_{i+1}}{1-x_i/x_{i+1}}(s_i-1)
\quad(i=1,\ldots ,n-1), 
\tag{1.5}
$$
where $s_i=(i,i+1)$ stands for the reflection with respect to the simple root
$\alpha_i$ $=\ep_i-\ep_{i+1}$. 
Note that 
$$
(T_i-t^{\half})(T_i+t^{-\half})=0\quad (i=1,\ldots ,n-1), 
\tag{1.6}
$$
and that 
the operators $T_1,\ldots ,T_{n-1}$ satisfy the fundamental 
relations for the canonical generators of the Hecke algebra $H(\frak S_n)$. 
Furthermore we define the $q$-{\it Dunkl operators} $Y_1,\ldots,Y_n$,
due to Cherednik,  by
$$
Y_i=T_iT_{i+1}\cdots T_{n-1}\omega T_1^{-1}\cdots T_{i-1}^{-1}
\quad(i=1,\ldots,n),
\tag{1.7}
$$
where 
$$
\omega=s_{n-1}\cdots s_1 \tau_1=\ldots =\tau_n s_{n-1}\cdots s_1. 
\tag{1.8}
$$
These operators $Y_1,\ldots ,Y_n$ commute with each other
and, for any {\it symmetric} Laurent polynomial 
$f(Y)$ of $Y=(Y_1,\ldots ,Y_n)$, 
Macdonald's symmetric polynomials $P_{\ld}(x)=P_{\ld}(x|q,t)$ 
$(\ld\in P^+)$ \cite{Ma2} satisfy the equation 
$$
f(Y)P_{\ld}(x) = P_{\ld}(x)f(q^\ld t^\rho),
\tag{1.9}
$$
where $f(q^\ld t^\rho)$ denotes the evaluation of $f$ at the point 
$q^\ld t^\rho=(q^{\ld_1} t^{\rho_1},\ldots ,q^{\ld_n} t^{\rho_n})$, 
and 
$$
\rho=\half\sum_{\alpha\in\Dt^+} \alpha=\half\sum_{i=1}^n (n-2i+1)\ep_i. 
\tag{1.10}
$$
\par\medpagebreak
One important fact is that the $q$-Dunkl operators have the triangularity with 
respect to a certain partial ordering of monomials. 
We define the partial ordering $\preceq$ in $P$ 
as follows: For $\ld,\mu\in P$, 
$$
\mu\preceq\ld\quad\Leftrightarrow\quad
\mu^+<\ld^+ \ \ \text{or}\ \ (\mu^+=\ld^+, \ \ \mu\le\ld),
\tag{1.11}
$$
where $\ld^+$ stands for the unique dominant integral weight in the $W$-orbit 
$W.\ld$ of $\ld$ and $\le$ is the dominance order 
($\mu\le\ld$ means that 
$\ld-\mu$ is a linear combination of positive roots 
with coefficients in $\N$). 
Then it turns out that,
for any $\ld,\mu\in P$, one has 
$$
Y^\mu x^\ld =x^\ld q^\br{\mu,\ld}t^\br{\mu,\rho(\ld)}+
\text{(lower order terms under $\preceq$)},
\tag{1.12}
$$
where
$$
\rho(\ld)=\half\sum_{\alpha\in \Delta^+}\vep(\br{\alpha,\ld})\alpha ,
\tag{1.13}
$$
where $\vep(u)=1$ if $u\ge0$ and $\vep(u)=-1$ if $u<0$.
Note that $\rho(\ld)$ is precisely the element $w_\ld.\rho$  in the $W$-orbit
of $\rho$ if one take the $w_\ld$ which has the minimal length among all 
$w\in W$ such that $\ld=w.\ld^+$. 
\Remark{1.1}{
Because of the definition of $q$-Dunkl operators mentioned above, 
the partial ordering $\preceq$ and the function $\vep(u)$ is different 
from those in \cite{Ma2}. 
Note that, under our definition of $\preceq$, 
the dominant weight $\ld^+$ is maximal in the $W$-orbit $W.\ld$. 
}
\medpagebreak
By the triangularity of $q$-Dunkl operators mentioned above, 
one can show that, for each $\ld\in P$, there exists a unique 
Laurent polynomial 
$E_\ld(x)=E_\ld(x|q,t)$ such that 
$$
E_\ld(x)=x^\ld+(\text{lower order terms under $\preceq$}), 
\tag{1.14}
$$
and that
$$
Y^\mu E_{\ld}(x)=E_\ld(x)q^\br{\mu,\ld} t^\br{\mu,\rho(\ld)}
\tag{1.15}
$$
for any $\mu\in P$. 
These Laurent polynomials $E_\ld(x)=E_\ld(x|q,t)$ are called the 
{\it nonsymmetric Macdonald polynomials} of type $A_{n-1}$. 
The first property (1.14) implies in particular 
that $E_\ld(x)$ is homogeneous of degree $|\ld|=\sum_{i=1}^n \ld_i$, 
and is eventually a polynomial in $x$ if $\ld\in L$. 
Note also that the second property (1.15) is equivalent to saying that 
$$
f(Y)E_{\ld}(x)=E_\ld(x)f(q^{\ld} t^{\rho(\ld)})
\tag{1.16}
$$
for {\it any} Laurent polynomial $f(Y)$ of the $q$-Dunkl operators. 
It is easy to see that the nonsymmetric Macdonald polynomials $E_\ld(x)$ 
($\ld\in P$) actually have coefficients in $\Q(q,t)$. 
We also remark that, as a function of $t$, 
each $E_\ld(x)=E_\ld(x|q,t)$ is regular 
at $t=q^k$ ($k=0,1,2,\ldots$). 
These polynomials $E_\ld(x)$ form a $\K$-basis of the ring $\K[x^{\pm1}]$ of 
Laurent polynomials or of the ring $\K[x]$ of polynomials as follows:
$$
\K[x^{\pm1}]=\bigoplus_{\ld\in P} \K E_\ld(x),\quad
\K[x]=\bigoplus_{\ld\in L} \K E_\ld(x). 
\tag{1.17}
$$
\par
It is known by \cite{Ma1} that, for any dominant integral weight $\ld\in P^+$, 
Macdonald's symmetric polynomial $P_\ld(x)$ is 
expressed as a linear combination 
of nonsymmetric Macdonald polynomials $E_\mu(x)$ ($\mu\in W.\ld$).  
To be more explicit, one has 
$$
P_\ld(x)=\sum_{\mu\in W.\ld} a_{\ld\mu} E_\mu(x) \quad\text{with}\quad
a_{\ld\mu}=\prod\Sb\alpha\in\Dt^+\\ \br{\alpha,\mu}<0 \endSb 
\frac{1-q^\br{\alpha,\mu}t^{\br{\alpha,\rho(\mu)}-1}}
{1-q^\br{\alpha,\mu}t^\br{\alpha,\rho(\mu)}}.
\tag{1.18}
$$
We also give a remark on the action of the Hecke algebra 
$H(\frak S_n)$ on nonsymmetric Macdonald polynomials: 
For each $i=1,\ldots n-1$,  one has 
$$
T_i E_\mu(x)=t^\half E_\mu(x) \quad \text{if}\quad \br{\alpha_i,\mu}=0,
\tag{1.19}
$$
and 
$$
T_i E_\mu(x)=u_{i,\mu} E_\mu(x)
+v_{i,\mu}E_{s_i\mu}(x) 
\quad\text{if}\quad \br{\alpha_i,\mu}\ne0.  
\tag{1.20}
$$
The coefficients $u_{i,\mu}$, $v_{i,\mu}$ are given by
$$
u_{i,\mu}=\frac{t^\half-t^{-\half}}
{1-q^{-\br{\alpha_i,\mu}}t^{-\br{\alpha_i,\rho(\mu)}}},\quad
v_{i,\mu}=t^{\half}
\tag{1.21}
$$
if $\br{\alpha_i,\mu}<0$, and
$$
\align
u_{i,\mu}=&\frac{t^\half-t^{-\half}}
{1-q^{-\br{\alpha_i,\mu}}t^{-\br{\alpha_i,\rho(\mu)}}}, \tag{1.22}\\
v_{i,\mu}=&
t^{-\half} 
\frac{(1-q^{\br{\alpha_i,\mu}}t^{\br{\alpha_i,\rho(\mu)}+1})
(1-q^{\br{\alpha_i,\mu}}t^{\br{\alpha_i,\rho(\mu)}-1})}
{(1-q^{\br{\alpha_i,\mu}}t^{\br{\alpha_i,\rho(\mu)}})^2}
\endalign
$$
if $\br{\alpha_i,\mu}>0$.  
\section{\S2: Formula of Cauchy type}
It is well-known that the Macdonald polynomials $P_\ld(x|q,t)$ ($\ld\in L^+$) 
have the following formula of Cauchy type \cite{Ma2}.  
Let now $x=(x_1,\ldots ,x_n)$ and $y=(y_1,\ldots ,y_n)$ be two sets of 
variables, and define the function $\Pi(x;y|q,t)$ by
$$
\Pi(x;y|q,t)=\prod_{1\le i,j\le n} 
\frac{(tx_iy_j;q)_\infty}{(x_iy_j;q)_\infty},
\tag{2.1}
$$
where $(x;q)_\infty=\prod_{i=0}^\infty(1-q^i x)$. 
The infinite products may be understood either 
in the sense of formal power series in appropriate variables, or 
in the analytic sense by assuming that 
$q$ is a complex number with $0<|q|<1$.  
Then we have 
$$
\Pi(x;y|q,t)=\sum_{\ld\in L^+} b_\ld(q,t) P_\ld(x|q,t) P_\ld(y|q,t),
\tag{2.2}
$$
where the coefficients are given by
$$
b_\ld(q,t)=\prod_{s\in\ld}\frac{1-q^{a(s)}t^{l(s)+1}}{1-q^{a(s)+1}t^{l(s)}}
\quad (\ld\in L^+)
\tag{2.3}
$$
in terms of the arm-length $a(s)=\ld_i-j$ and 
the leg-length $l(s)=\ld'_j-i$ of 
each box $s=(i,j)$ in the partition $\ld$. 

\par\medpagebreak

We now introduce the function $E(x;y|q,t)$ by setting 
$$
E(x;y|q,t)=
\prod_{1\le j<i\le n}\frac{(qtx_iy_j;q)_\infty}{(qx_iy_j;q)_\infty}
\prod_{1\le i\le n}   \frac{(qtx_iy_i;q)_\infty}{(x_iy_i;q)_\infty}
\prod_{1\le i<j\le n}\frac{(tx_iy_j;q)_\infty}{(x_iy_j;q)_\infty}.
\tag{2.4}
$$
Note that this function can be factored as follows: 
$$
E(x;y|q,t)=\Pi(x;y|q,t) 
\prod_{i=1}^n\frac{1}{1-tx_iy_i}
\prod_{j<i}\frac{1-x_iy_j}{1-tx_iy_j}. 
\tag{2.5}
$$
The ratio $E(x;y|q,t)\Pi(x;y|q,t)^{-1}$ is essentially one of 
the rational functions 
(before symmetrization) which are used in \cite{Mi1} and \cite{KN}. 
\Theorem{2.1}{
The function $E(x;y|q,t)$ has an expansion 
$$
E(x;y|q,t)=\sum_{\ld\in L} a_{\ld}(q,t) E_\ld(x|q,t) E_\ld(y|q^{-1},t^{-1})
\quad(a_{\ld}(q,t)\in \Q(q,t))
\tag{2.6}
$$
summed over all compositions $\ld\in L$.
}
\par
In order to describe the coefficients $a_\ld(q,t)$ 
$(\ld\in L)$ in the expansion (2.6), 
we use the notion of leg-length generalized to compositions, 
due to Knop and Sahi \cite{KS}. 
For each box $s=(i,j)$ in a composition $\ld\in L$, the {\it generalized 
leg-length} $l(s)=l_\ld(s)$ is defined to be the sum 
$$
l(s)=l_{\text{up}}(s)+l_{\text{low}}(s)
\tag{2.7}
$$
of the upper and the lower leg-length, where
$$
l_{\text{low}}(s)=\#\{k>i~;~j\le\ld_k\le\ld_i\},\quad
l_{\text{up}}(s)=\#\{k<i~;~j\le\ld_k+1\le\ld_i\}.
\tag{2.8}
$$
Note that this $l(s)$ is the same as the usual 
leg-length if $\ld$ is a partition. 
\Theorem{2.2}{
For each composition $\ld\in L$, 
the coefficient $a_\ld(q,t)$ in expansion $(2.6)$ is given by
$$
a_\ld(q,t)=\prod_{s\in\ld}\frac{1-q^{a(s)+1}t^{l(s)+1}}{1-q^{a(s)+1}t^{l(s)}}
\quad (\ld\in L), 
\tag{2.9}
$$
where $l(s)=l_\ld(s)$ $(s\in\ld)$ stands for the 
generalized leg-length in $\ld$.
}
\noindent
In this section, we will give a proof of Theorem 2.1.  
Also, we will describe the ratio of $a_\ld(q,t)$ and $a_\mu(q,t)$ 
when $\ld$ is a partition and $\mu$ is a composition in the orbit $W.\ld$.  
Theorem 2.2 will be established in two ways in Sections 3 and 4, 
by determining $a_\ld(q,t)$ for partitions $\ld\in L^+$.  

\par\medpagebreak

In what follows, we set $K(x;y|q,t)=E(x;y^{-1}|q,t)$, i.e. 
$$
K(x;y|q,t)=
\prod_{1\le j<i\le n}\frac{(qtx_i/y_j;q)_\infty}{(qx_i/y_j;q)_\infty}
\prod_{1\le i\le n}   \frac{(qtx_i/y_i;q)_\infty}{(x_i/y_i;q)_\infty}
\prod_{1\le i<j\le n}\frac{(tx_i/y_j;q)_\infty}{(x_i/y_j;q)_\infty}.
\tag{2.10}
$$
\Proposition{2.3}{
For each $i=1,\ldots ,n$, one has
$$
Y_{i,x} K(x;y|q,t) = (Y_{i,y}^\ast)^{-1} K(x;y|q,t),
\tag{2.11}
$$
where the suffix $x$ or $y$ indicates the variables 
on which the operator should 
act, and $Y_i^\ast$ is the dual $q$-Dunkl operator 
$($ cf. \cite{KN}$)$ defined by
$$
Y_i^\ast=T_i^{-1}\cdots T_{n-1}^{-1}\omega T_1\cdots T_{i-1}.
\tag{2.12}
$$
}
\noindent
This proposition is a direct consequence of the following lemma. 
\Lemma{2.4}{\newline
$(1)$ \quad $T_{i,x} K(x;y|q,t)=T_{i,y} K(x;y|q,t)$
\quad$(i=1,\ldots ,n-1)$. \newline
$(2)$ \quad $\omega_x K(x;y|q,t)=\omega_y^{-1} K(x;y|q,t)$.
}
\noindent
Note that the function $K(x;y|q,t)$ can be factored as follows: 
$$
K(x;y|q,t)=\Pi(x;y^{-1}|q,t)\psi(x,y),\ \ 
\psi(x,y)=
\prod_{i=1}^n\left(\frac{1}{1-tx_i/y_i}
\prod_{j<i}\frac{1-x_i/y_j}{1-tx_i/y_j}\right). 
\tag{2.13}
$$
Since $\Pi(x;y^{-1}|q,t)$ is symmetric both in $x$ and in $y$, 
the formula of Lemma 2.4.(1) is equivalent to 
$$
T_{i,x}\psi(x,y)=T_{i,y}\psi(x,y)\quad(i=1,\ldots ,n-1). 
\tag{2.14}
$$
For a fixed $i$, it reduces to the identity
$$
T_{i,x}\psi_{i,i+1}(x,y)=T_{i,y}\psi_{i,i+1}(x,y)
\tag{2.15}
$$
for 
$$
\psi_{i,i+1}(x,y)=
\frac{1-x_{i+1}/y_i}{(1-tx_i/y_i)(1-tx_{i+1}/y_i)(1-tx_{i+1}/y_{i+1})},
\tag{2.16}
$$
which can be checked by a direct computation. 
(This computation is essentially 
contained in \cite{Mi1}, \cite{MN}).  
The formula of Lemma 2.4.(2) can be proved directly 
by chasing the arguments of 
$q$-shift factorials under the action $\omega_y\omega_x$. 

\par

\Proofof{Theorem 2.1}{
We begin with expanding the function $E(x;y|q,t)$ in the form
$$
E(x;y|q,t)=\sum_{\ld\in L}E_\ld(x|q,t) F_\ld(y|q,t)\quad
(F_\ld(y|q,t)\in \Q(q,t)[y]).
\tag{2.17}
$$
We will show that each $F_\ld(x|q,t)$ is a constant multiple of 
$E_\ld(y|q^{-1},t^{-1})$. 
Note that 
Proposition 2.3 implies 
$$
Y_x^{-\mu} K(x;y|q,t)=(Y^\ast_y)^{\mu}K(x;y|q,t)
\tag{2.18}
$$ 
for any $\mu\in P$. 
Since 
$$
K(x;y|q,t)=\sum_{\ld\in L}E_\ld(x|q,t) F_\ld(y^{-1}|q,t),
\tag{2.19}
$$
we have 
$$
{Y^\ast_y}^{\mu} F_\ld(y^{-1}|q,t)= 
F_\ld(y^{-1}|q,t) q^{-\br{\mu,\ld}}t^{-\br{\mu,\rho(\ld)}}. 
\tag{2.20}
$$
As is shown in \cite{KN}, for each $i=1,\ldots ,n$ the $q$-Dunkl operator
$Y_i$ and its dual $Y^\ast_i$ are interchanged by the involution 
$\iota$ on $\K[y]$
such that $\iota(y_j)=y_j^{-1}$ $(j=1,\ldots ,n)$, $\iota(q)=q^{-1}$, 
$\iota(t^{\half})=t^{-\half}$. 
Hence we have 
$$
Y_y^{\mu} F_\ld(y|q^{-1},t^{-1})= 
F_\ld(y|q^{-1},t^{-1}) q^{\br{\mu,\ld}}t^{\br{\mu,\rho(\ld)}}
\tag{2.21}
$$
for  all $\mu\in P$. 
This implies that $F_\ld(y|q^{-1},t^{-1})$ is a constant multiple of 
$E_\ld(y|q,t)$.  
Namely, $F_\ld(y|q,t)$ is a constant multiple of 
$E_\ld(y|q^{-1},t^{-1})$.  
}
\par\medpagebreak
Before determining the coefficients $a_\ld(q,t)$, 
we will describe the relation 
between $a_\ld(q,t)$ and $a_\mu(q,t)$ when $\ld$ is dominant and $\mu$ 
is in the orbit $W.\ld$.  
\Lemma{2.5}{
If $\ld\in L^+$, $\mu\in L$ and $\mu\in W.\ld$, then one has
$$
a_\mu(q,t)=a_\ld(q,t)\prod\Sb \alpha\in\Dt^+\\ \br{\alpha,\mu}<0\endSb
\frac{(1-q^{-\br{\alpha,\mu}}t^{-\br{\alpha,\rho(\mu)}+1})
(1-q^{-\br{\alpha,\mu}}t^{-\br{\alpha,\rho(\mu)}-1})}
{(1-q^{-\br{\alpha,\mu}}t^{-\br{\alpha,\rho(\mu)}})^2}
\tag{2.22}
$$
}
\Proof{
By Lemma 2.4.(1), we have
$$
\sum_{\mu\in L}a_\mu T_{i,x}E_\mu(x) \iota(E_\mu(y))
=\sum_{\mu\in L}a_\mu E_\mu(x) T_{i,y} \iota(E_\mu(y)),
\tag{2.23}
$$
for each $i=1,\ldots ,n-1$, where $a_\mu=a_\mu(q,t)$. 
As we remarked at the end of Section 2, 
for each $\mu\in L$, we have 
$
T_{i,x}E_\mu(x)=t^\half E_\mu(x)$ if $\br{\alpha_i,\mu}=0$, 
and 
$$
T_{i,x}E_\mu(x)= u_{i,\mu} E_{\mu}(x) + v_{i,\mu} E_{s_i\mu}(x).  
\tag{2.24}
$$
when $\br{\alpha_i,\mu}\ne 0$.  
Since 
$T_{i,y}\iota(E_\mu(y))=\iota(T_{i,y}^{-1}E_\mu(y))$, 
and $T_{i,y}^{-1}=T_{i,y}-(t^\half-t^{-\half})$,
we can determine the action of $T_{i,y}$ on $\iota(E_\mu(y))$ as
follows: 
$T_{i,\mu}\iota(E_\mu(y))=t^{\half}\iota(E_\mu(y))$
if $\br{\alpha_i,\mu}=0$, and 
$$
T_{i,y}\iota(E_\mu(y))
=(\iota(u_{i,\mu})+(t^\half-t^{-\half}) )\iota(E_\mu(y))
+\iota(v_{i,\mu})\iota(E_{s_i\mu}(y))
\tag{2.25}
$$
if $\br{\alpha_i,\mu}\ne0$.
By substituting these formulas into (2.23), we obtain 
the recurrence formulas
$$
a_{s_i\mu} v_{s_i\mu}=a_\mu \iota(v_\mu) 
\tag{2.26}
$$
for $\mu\in L$ with $\br{\alpha_i,\mu}\ne0$.   
Hence, by  (2.21) and (2.22), we have
$$
a_\mu=a_{s_i\mu}
\frac{(1-q^{-\br{\alpha_i,\mu}}t^{-\br{\alpha_i,\rho(\mu)}+1})
(1-q^{-\br{\alpha_i,\mu}}t^{-\br{\alpha_i,\rho(\mu)}-1})}
{(1-q^{-\br{\alpha_i,\mu}}t^{-\br{\alpha_i,\rho(\mu)}})^2}
\tag{2.27}
$$
for $\mu\in L$ with $\br{\alpha_i,\mu}<0$. 
Assume now that $\ld\in L^+$ is a partition and $\mu\in W.\ld$.  
Then one can find a sequence of 
simple roots $\alpha_{j_1},\ldots ,\alpha_{j_p}$
such that
$\mu=s_{j_1}\cdots s_{j_p}(\ld)$ and that
$$
\br{\alpha_{j_1},\mu}<0,
\br{\alpha_{j_2},s_{j_1}(\mu)}<0,
\ldots ,
\br{\alpha_{j_p},s_{j_{p-1}}\cdots s_{j_1}(\mu)}<0. 
\tag{2.28}
$$
Note also that 
$$
\{\alpha\in\Dt^+~;~\br{\alpha,\mu}<0\}
=\{\alpha_{j_1},s_{j_1}(\alpha_{j_2}),\ldots,
s_{j_1}\cdots s_{j_{p-1}}(\alpha_{j_p})\}.  
\tag{2.29}
$$
Applying formula (2.27) to $\mu^{(0)}=\mu, \mu^{(1)}=s_{j_1}(\mu),\cdots, 
\mu^{(p)}=s_{j_p}\cdots s_{j_1}(\mu)=\ld$ repeatedly, we obtain 
formula (2.22) since 
$\br{\alpha_{j_r},\mu^{(r-1)}}=
\br{s_{j_1}\cdots s_{j_{r-1}}(\alpha_{j_r}),\mu}$
and  
$\br{\alpha_{j_r},\rho(\mu^{(r-1)})}=
\br{s_{j_1}\cdots s_{j_{r-1}}(\alpha_{j_r}),\rho(\mu)}$
for $r=1,\ldots ,p$.
}
Lemma 2.5 can be rewritten in the combinatorial language. 
Imitating Sahi's notation \cite{S}, we set
$$
d_\mu(q,t)=\prod_{s\in\mu}(1-q^{a(s)+1}t^{l(s)+1}),\quad
d'_\mu(q,t)=\prod_{s\in\mu}(1-q^{a(s)+1}t^{l(s)})
\tag{2.30}
$$
for each $\mu\in L$.
In this notation, 
Theorem 2.2 is equivalent to the equality $a_\mu(q,t)=d_\mu(q,t)/d'_\mu(q,t).$
\Lemma{2.6}{
If $\ld\in L^+$, $\mu\in L$ and $\mu\in W.\ld$, then formula $(2.22)$ has 
an alternative expression 
$$
a_\mu(q,t)=a_\ld(q,t) \frac{d'_\ld(q,t)d_\mu(q,t)}{d_\ld(q,t)d'_\mu(q,t)}. 
\tag{2.31}
$$
}
\Proof{
For each $\mu\in L$, set $a'_\mu=a'_\mu(q,t)=d_\mu(q,t)/d'_\mu(q,t)$.  
For the proof of formula (2.31), it is enough to show
$$
a_\mu(q,t)=a_{s_i\mu}(q,t) \frac{a'_\mu(q,t)}{a'_{s_i\mu}(q,t)} 
\tag{2.32}
$$
assuming that $\br{\alpha_i,\mu}<0$ ($i=1,\ldots ,n-1$); 
one can use (2.32) repeatedly to prove (2.31) in view of 
the expression $\mu=s_{j_1}\cdots s_{j_p}(\ld)$ as in the 
proof of Lemma 2.5. 
When $\br{\alpha_i,\mu}<0$, it is easy to see that the only difference 
between $a'_\mu(q,t)$ and $a'_{s_i\mu}(q,t)$
arises at the box $s=(i+1,\mu_i+1)\in\mu$ (or at $s'=(i,\mu_i+1)\in s_i\mu$). 
In fact we have
$$
a'_\mu=a'_{s_i\mu}
\frac{(1-q^{\mu_{i+1}-\mu_i}t^{l_\mu(s)+1})
(1-q^{\mu_{i+1}-\mu_i}t^{l_\mu(s)-1})}
{(1-q^{\mu_{i+1}-\mu_i}t^{l_\mu(s)})^2}.
\tag{2.33}
$$
On the other hand, one can directly check that
$$
-\br{\alpha_i,\mu}=\mu_{i+1}-\mu_i,\quad
-\br{\alpha_i,\rho(\mu)}=l_\mu(s)
\tag{2.34}
$$
by the definition (1.13) of $\rho(\mu)$.  
Hence we have (2.32) by comparing (2.27) and (2.33).
}

\section{\S3: First proof of Theorem 2.2}
In this section, we calculate the coefficients $a_\ld=a_\lambda(q,t)$ for
partitions $\lambda\in L^{+}$ by means of 
asymptotic analysis of a $q$-Selberg type integral.   
Such an argument has been employed in \cite{Mi2} to evaluate the
scalar product for Macdonald's symmetric polynomials. %
\par\medpagebreak
We now assume that $q$ and $t$ are complex numbers with $0<|q|,|t|<1$, 
and recall Cherednik's scalar product \cite{C1}.
For $f=f(x|q,t), g=g(x|q,t)\in {\Bbb K}[x]=
{\Bbb K}[x_1,\ldots,x_n]$, we define
$$
(f,g) =\left(\frac{1}{2\pi\sqrt{-1}}\right)^n
\int_{\Bbb T^n}f(x|q,t)g(x^{-1}|q^{-1},t^{-1})w(x|q,t)
\frac{dx_1\cdots dx_n}{x_1\cdots x_n}\,,
\tag{3.1}
$$ 
where
$$
w(x|q,t)=\prod_{1\le i<j\le n}
\frac{(x_i/x_j;q)_\infty(qx_j/x_i;q)_\infty}
{(tx_i/x_j;q)_\infty(qtx_j/x_i;q)_\infty}
\tag{3.2}
$$
and
$\Bbb T^n=\{(x_1,\ldots,x_n)\in {\Bbb C}^n ; |x_i|=1 \ (1\le i\le n)\}$
with the standard orientation.  
Note that
$$
w(x|q,q^k)=\prod_{1\le i<j\le n}(x_i/x_j;q)_k(qx_j/x_i;q)_k 
\tag{3.3} 
$$
if $t=q^k$ ($k=0,1,2,\ldots$), where 
$(a;q)_k=(a;q)_\infty/(at;q)_\infty=\prod_{0\le i\le k-1}(1-aq^i)$.
The nonsymmetric Macdonald polynomials $E_\ld(x;q,t)$ 
($\ld\in L$) form an orthogonal basis of $\K[x]$ with respect 
to this scalar product: 
$$
(E_\lambda, E_\mu)=0\quad \text{if}\quad
\ld\ne\mu\,.
\tag{3.4}
$$  
It is known furthermore that, if $t=q^k$ ($k\in \N$) and $\ld\in L^+$ is 
a partition, then 
$$
(E_\ld , E_\ld)=\prod_{1\le i< j\le n }
\frac{(q^{\ld_i-\ld_{j}+1+k(j-i)};q)_k}
{(q^{\ld_i-\ld_{j}+1+k(j-i-1)};q)_{k}}.
\tag{3.5}
$$ 
(See \cite{Ma1}, \cite{C1}.) 
For general values of $t$ with $|t|<1$, one has
$$
(E_\ld , E_\ld)=\prod_{1\le i< j\le n }
\frac{(q^{\ld_i-\ld_{j}+1}t^{j-i};q)_\infty^2}
{(q^{\ld_i-\ld_{j}+1}t^{j-i+1};q)_{\infty}
(q^{\ld_i-\ld_{j}+1}t^{j-i-1};q)_{\infty}}
\tag{3.6}
$$ 
for any partition $\ld\in L^+\,$.
Note that, as functions in $t$, the both sides of (3.6) are 
meromorphic in $\{ |t|<1 \}$.  
Since this formula is valid at the points $t=q^k$ ($k\in\N$) accumulating 
at the origin, one can conclude that the left hand side of (3.6) 
is eventually holomorphic 
near $t=0$, and that (3.6) is valid as an equality of analytic functions. 
It is also known that, if $\mu\in W.\ld$ is a composition in the $W$-orbit of 
a partition $\ld$, then we have
$$
\frac{(E_\mu, E_\mu)}{(E_\ld,E_\ld)}=
\prod\Sb\alpha\in\Dt^+\\ \br{\alpha,\mu}<0\endSb
\frac{(1-q^{-\br{\alpha,\mu}}t^{-\br{\alpha,\rho(\mu)}})^2}
{(1-q^{-\br{\alpha,\mu}}t^{-\br{\alpha,\rho(\mu)}+1})
(1-q^{-\br{\alpha,\mu}}t^{-\br{\alpha,\rho(\mu)}-1})}. 
\tag{3.7}
$$
\par\medpagebreak
We now consider $x=(x_1,\ldots,x_n)$  as variables inside the polydisc 
$\{|x_j|<1~;~j=1,\ldots,n \}\subset\C^n$.  
Note that the series expansion 
$$
K(x;y|q,t)=\sum_{\mu\in L}
a_{\mu}(q,t)E_{\mu}(x|q,t)E_{\mu}(y^{-1}|q^{-1},t^{-1})\,,
\tag{3.8}
$$ 
in Theorem 2.1 is then uniformly convergent on $\Bbb T^n$ in $y$.  
Hence by the orthogonality relations (3.4) we have 
$$\align
&\left(\frac{1}{2\pi\sqrt{-1}}\right)^n
\int_{\Bbb T^n} K(x;y|q,t)E_\ld(y|q,t)w(y|q,t)
\frac{dy_1\cdots dy_n}{y_1\cdots y_n}\tag{3.9}\\
&=\sum_{\mu\in L} a_{\mu}(q,t)E_{\mu}(x|q,t)
(E_{\ld}, E_{\mu})\\
&= a_{\ld}(q,t) E_{\ld}(x|q,t)
(E_{\ld}, E_{\ld})\,.
\endalign$$
We study the asymptotic behavior of the left hand side of (3.9) in $x$ 
in the region
$$
1>|x_1|\gg |x_2|\gg\cdots\gg |x_n|\,.\tag{3.10} 
$$
For the moment, we assume that $t=q^k$ ($k=0,1,2,\ldots$) and that 
$\ld$ is a partition. 
\Proposition{3.1}{
If $t=q^k$ ($k\in\N$) and $\ld\in L^+\,,$ then one has
$$
\align
&\left(\frac{1}{2\pi\sqrt{-1}}\right)^n
\int_{{\Bbb T}^n} K(x;y|q,q^k)E_\ld(y|q,q^k)w(y|q,q^k)
\frac{dy_1\cdots dy_n}{y_1\cdots y_n}
\tag{3.11}\\
&\quad\sim\quad
x^{\ld}~
\prod_{i=1}^{n}\frac{(q^{\ld_i+(n-i)k+1};q)_k}{(q;q)_{k}},
\endalign
$$
as $\max\{|x_2/x_1|,|x_3/x_2|,\ldots,|x_n/x_{n-1}|\}\,$
tends to $0$. 
}
\noindent
Proposition will be proved later in this section.
\par\medpagebreak
We apply Proposition 3.1 to compare the coefficients of $x^\ld$ 
in equality (3.9).  
Since we know by (3.9) that the integral has the asymptotic behavior 
$x^\ld a_\ld(q,t)(E_\ld,E_\ld)+\cdots$, 
we obtain 
$$
a_{\ld}(q,q^k)\,(E_{\ld}, E_{\ld})
=\prod_{i=1}^{n}\frac{(q^{\ld_i+(n-i)k+1};q)_k}{(q;q)_{k}}
\tag{3.12}
$$
for any partition $\ld\in L^+$, provided that $t=q^k$ ($k\in\N$). 
This is equivalent to
$$
\align
a_{\ld}(q,q^k)
&=\prod_{1\le i\le j\le n }
\frac{(q^{\ld_i-\ld_{j+1}+1+k(j-i)};q)_k}
{(q^{\ld_i-\ld_{j}+1+k(j-i)};q)_{k}}
\tag{3.13}\\
&=\prod_{i=1}^n\prod_{j=i }^n
\frac{(q^{\ld_i-\ld_{j}+1+k(j-i+1)};q)_{\ld_j-\ld_{j+1}}}
{(q^{\ld_i-\ld_{j}+1+k(j-i)};q)_{\ld_j-\ld_{j+1}}}
\endalign
$$
by formula (3.5), where  we set $\ld_{n+1}=0$. 
By analytic continuation as before, we get
$$
a_{\ld}(q,t)
=\prod_{i=1}^n\prod_{j=i }^n
\frac{(q^{\ld_i-\ld_{j}+1}\,t^{j-i+1};q)_{\ld_j-\ld_{j+1}}}
{(q^{\ld_i-\ld_{j}+1}\,t^{j-i};q)_{\ld_j-\ld_{j+1}}},
\tag{3.14}
$$
since $a_\ld(q,t)$ is a rational function in $t$. 
Formula (3.14) implies that $a_0=1$ and that
$$
a_{\ld+(1^m)}(q,t)=a_\ld(q,t)
\prod_{i=1}^m
\frac{1-q^{\ld_i+1}t^{m-i+1}}{1-q^{\ld_i+1}t^{m-i}}
\tag{3.15}
$$
for any $\ld\in L^+$ with $l(\ld)\le m$.
In fact, the difference between $a_\ld(q,t)$ and $a_{\ld+(1)^m}(q,t)$ 
appears only in the factors in (3.14) with $1\le i\le m$ and $j=m$. 
{}From (3.15) it follows immediately that
$$
a_\ld(q,t)=\prod_{s\in\ld}\frac{1-q^{a(s)+1}t^{l(s)+1}}{1-q^{a(s)+1}t^{l(s)}}
=\frac{d_\ld(q,t)}{d'_\ld(q,t)}
\tag{3.16}
$$
for all $\ld\in L^+$.
Hence by Lemma 2.6, the same formula (3.16) holds for all compositions 
$\ld\in L$ if $l(s)$ is understood as the generalized leg-length.  
This completes the proof of Theorem 2.2. 
\par\medpagebreak
Note that formula (3.12) extends to the equality
$$
a_{\ld}(q,t)\,(E_{\ld}, E_{\ld})
=\prod_{i=1}^{n}\frac{(q^{\ld_i+1}t^{n-i};q)_\infty (qt;q)_\infty}
{(q^{\ld_i+1}t^{n-i+1};q)_\infty(q;q)_\infty}\quad(\ld\in L^+)
\tag{3.17}
$$
of analytic functions in $t$. 
On the other hand, by comparing formula (2.22) of Lemma 2.5 and (3.7), 
we see that 
$$
a_{\mu}(q,t)\,(E_{\mu}, E_{\mu})=a_{\ld}(q,t)\,(E_{\ld}, E_{\ld})
\tag{3.18}
$$
for all compositions $\mu\in W.\ld$.  
Summarizing these remarks, we obtain the following theorem.
\Theorem{3.2}{
For each composition $\ld\in \N^n$, we have 
$$
\left(\frac{1}{2\pi\sqrt{-1}}\right)^n\int_{\Bbb T^n}
 K(x;y|q,t) E_\ld(y|q,t)w(y|q,t)
\frac{dy_1\cdots dy_n}{y_1\cdots y_n}=C_\ld. E_\ld(x|q,t)
\tag{3.19}
$$
for $x=(x_1,\ldots,x_n)\in \C^n$ with $|x_j|<1$ $(j=1,\ldots,n)$.
Here $C_\ld$ is constant on each $W$-orbit, and is given by 
$$
C_\ld=\left(\frac{(qt;q)_\infty}{(q;q)_\infty}\right)^n
\prod_{i=1}^n \frac{(q^{\ld_i+1}t^{n-i};q)_\infty}{(q^{\ld_i+1}t^{n-i+1};q)_\infty}
\tag{3.20}
$$
when $\ld\in L^+$ is a partition.
}
\noindent
In this sense, our function $K(x;y|q,t)$ is a reproducing kernel 
for nonsymmetric Macdonald polynomials. 
\par\medpagebreak
In the rest of this section, we will prove Proposition 3.1. 
{}From now on we assume that $t=q^k$ for a fixed $k\in\N$, and omit 
$(q,t)=(q,q^k)$ in the notation unless it might lead to confusion. 
\par
In proving Proposition, we may assume that $0<|x_j|<1$ 
for $j=1\ldots,n$ and that all $x_j$'s are mutually distinct. 
Recall that 
$$
\align
K(x;y)&=
\prod_{j<i}\frac{1}{(qx_i/y_j;q)_k}
\prod_{i}   \frac{1}{(x_i/y_i;q)_{k+1}}
\prod_{ i<j}\frac{1}{(x_i/y_j;q)_k}
\tag{3.21}\\
&=\prod\Sb  i,j \endSb \frac{1}{(q^{\theta(i>j)}x_i/y_j;q)_{k+\dt_{ij}}},
\endalign
$$
where $\theta(i>j)=1 $ if $i>j$, and $\theta(i>j)=0$ if $i\le j$, and
$$
w(y)=\prod_{i<j}(y_i/y_j;q)_k (qy_j/y_i)_k. 
\tag{3.22}
$$
Note first that, as a function of $y_j$ ($1\le j\le n$), the integrand 
$$
K(x;y)E_\ld(y)w(y)(y_1\cdots y_n)^{-1}
\tag{3.23}
$$
of (3.11) is regular at $y_j=0$ and has poles only at 
$y_j=x_sq^{l}$ ($1\le s\le n, l\in\N $). 
The range of $l$ can be specified as follows:  
$$
\matrix
\text{\rm (1)}&\quad&0\le l <k &\quad\text{if}\quad &1\le s<j, &\qquad\\ 
\text{\rm (2)}&\quad&0\le l \le k &\quad\text{if}\quad  &s=j,  &\quad\text{and} \\ 
\text{\rm (3)}&\quad&0<l\le k &\quad\text{if}\quad  & j<s\le n. &\qquad 
\endmatrix
\tag{3.24}
$$
The integral (3.11) will be computed by picking up successively 
the residues at $y_j=x_s q^l$ $(1\le j,s\le n)$ with $l$ satisfying 
(3.24). 
To make clear this inductive process, we propose a lemma. 
\par
For any pair  $(I,J)$ of subsets of $\{1,\ldots,n\}$ with $|I|=|J|$, 
we extend the notation of (3.21), (3.23) as follows: 
$$
K(x_I;y_J)=\prod\Sb  i\in I\\j\in J \endSb 
\frac{1}{(q^{\theta(i>j)}x_i/y_j;q)_{k+\dt_{ij}}},
\ \ 
w(y_J)=\prod\Sb i,j\in J\\ i<j\endSb
(y_i/y_j;q)_k (qy_j/y_i)_k,
\tag{3.25} 
$$
where $x_I=(x_i)_{i\in I}$ and $y_J=(y_j)_{j\in J}$.
\Lemma{3.3}{
Fix two indices $j\in J$, $s\in I$ and $l\in\N$ satisfying $(3.24)$, 
and set $J'=J\backslash \{j\}$ and $I'=I\backslash\{s\}$. 
Then the residue
$$
f(x_I;y_{J'})=
\operatorname{Res}_{y_j=x_sq^l}(K(x_I;y_J)w(y_J)\frac{dy_j}{y_j})
\tag{3.26}
$$
at $y_j=x_sq^l$ has an expression 
$$
f(x_I;y_{J'})=g(x_I;y_{J'}) K(x_{I'};y_{J'}) w(y_{J'}) , 
\tag{3.27}
$$
where $g(x_I;y_{J'})$ is a polynomial in $y_{J'}$ with coefficients 
in $\K(x_I)$. 
} 
\Proof{
The only factors to be checked are 
$$
\frac{(q^{-l}y_\mu/x_s;q)_k (q^{l+1}x_s/y_\mu;q)_k}
{(q^{\theta(s>\mu)}x_s/y_\mu;q)_{k+\dt_{s\mu}}}
\tag{3.28}
$$
for $\mu\in J'$ with $\mu<j$, and 
$$
\quad
\frac{(q^lx_s/y_\mu;q)_k (q^{1-l}y_\mu/x_s;q)_k}
{(q^{\theta(s>\mu)}x_s/y_\mu;q)_{k+\dt_{s\mu}}}
\tag{3.29} 
$$
for $\mu\in J'$ with $\mu>j$. 
As a function of $y_\mu$, the numerators of (3.28) and (3.29) have
zeros at $y_\mu=x_s q^{m}$ for $m\in\{l-k+1,\l-k+2,\cdots,l+k\}$ 
and for $m\in \{l-k,\l-k+1,\cdots,l+k-1\}$, respectively. 
{}From this, it turns out that the rational functions (3.28) 
and (3.29) are in fact regular at $y_\mu=x_s q^m$ for all $m\in\N$, 
provided that $l$ satisfies the condition of (3.24). 
}
Let us now integrate $K(x;y)E_\ld(y)w(y)/y_1$ with respect to the 
variable $y_1$.
Then we have 
$$
\align
&\frac{1}{2\pi\sqrt{-1}}\int_{|y_1|=1} 
K(x;y)E_\ld(y)w(y)\frac{dy_1}{y_1}
\tag{3.30}\\
&\quad=\sum_{s=1}^n \sum_{l=0}^k 
\operatorname{Res}_{y_1=x_sq^l}(K(x;y)E_\ld(y)w(y)\frac{dy_1}{y_1})
\endalign
$$
since the integrand is regular at $y_1=0$. 
If we regard each summand of the right-hand side as a function of 
$y_2$, it has a zero at $y_2=0$ and has poles only at 
$y_2=x_r q^m$ for $r\ne s$, $0\le m\le k$ by Lemma 3.3. 
Hence we have
$$
\align
&\left(\frac{1}{2\pi\sqrt{-1}}\right)^2 \int_{|y_1|=|y_2|=1} 
K(x;y)E_\ld(y)w(y)\frac{dy_1 dy_2}{y_1 y_2}
\tag{3.31}\\
&\quad=\sum\Sb 1\le s,r\le n \\ s\ne r \endSb 
\sum_{0\le l,m\le k} 
\operatorname{Res}\Sb y_1=x_sq^l \\ \ y_2=x_rq^m\endSb
(K(x;y)E_\ld(y)w(y)\frac{dy_1dy_2}{y_1y_2}).
\endalign
$$
Applying Lemma 3.3 repeatedly from $j=1$ to $n,$ 
we obtain the equality 
$$
\align
&\left(\frac{1}{2\pi\sqrt{-1}}\right)^n
\int_{\Bbb T^n} K(x;y)E_\ld(y)w(y)
\frac{dy}{y}\tag{3.32}\\
&=\sum_{\sigma\in{\frak S}_n}\sum_{0\le l_1,\cdots,l_n\le k}
\text{Res}_{y=(x_{\sigma(1)}q^{l_1},\ldots, x_{\sigma(n)}q^{l_n})}
\left(K(x;y)E_\ld(y)w(y)
\frac{dy}{y}\right),
\endalign
$$
where we used the abbreviation $dy/y=dy_1\cdots dy_n/y_1\cdots y_n$.
\par
We now investigate the asymptotic behavior of this 
function in the region (3.10). 
Note that we have
$$
\align
&\text{Res}_{y=(x_{\sigma(1)}q^{l_1},\ldots, x_{\sigma(n)}q^{l_n})}
\left(K(x;y)E_\ld(y)w(y)
\frac{dy}{y}\right)
\tag{3.33}\\
&=\text{Res}_{y=(q^{l_1},\ldots, q^{l_n})}
\left(K(x;x_{\sigma}y)E_\ld(x_{\sigma}y)w(x_{\sigma}y)
\frac{dy}{y}\right),
\endalign 
$$
where $x_{\sigma}y$ stands for 
$(x_{\sigma(1)}y_1,\ldots, x_{\sigma(n)}y_n).$
The function $K(x;x_\sigma y)w(x_\sigma y)$ 
can be written in the following form:
$$
\prod\Sb  i,j \endSb \frac{1}{(q^{\theta(i\ge j)}x_i/x_{\sigma(j)}y_j;q)_{k}}
\prod_{i\ne j} (q^{\theta(i>j)}x_{\sigma(i)}y_i/x_{\sigma(j)}y_j;q)_k 
\times \prod_{j=1}^n\frac{1}{1-x_j/x_{\sigma(j)}y_j}. 
\tag{3.34}
$$
The product of the first two factors altogether is 
bounded in the limit (3.10).  
If $\sigma\in\frak S_n$ is {\it not} the identity element, 
one can take a suffix $i$ such that $i<\sigma(i).$
As an effect of the factor $1/(1-x_i/x_{\sigma(i)}y_i)$ in the third 
factors of (3.34), we then have
$$
|K(x;x_\sigma y)w(x_\sigma y)|
=O\left(\left|\frac{x_{\sigma(i)}}{x_{i}}\right|\right). 
\tag{3.35}
$$
Since $\ld$ is a partition, $x^{-\ld}E_\ld(x_\sigma y)$ is also 
bounded in the limit (3.10). 
Hence we have 
$$
x^{-\ld}\text{Res}_{y=(q^{l_1},\ldots, q^{l_n})}
\left(K(x;x_{\sigma}y)E_\ld(x_{\sigma}y)w(x_{\sigma}y)
\frac{dy}{y}\right)=
O\left(\left|\frac{x_{\sigma(i)}}{x_{i}}\right|\right), 
\tag{3.36}
$$
provided that $\sigma$ is not the identity element. 
If $\sigma$ is the identity element, the function
$$
\text{Res}_{y=(q^{l_1},\ldots, q^{l_n})}
\left(K(x;x y)E_\ld(x y)w(x y)
\frac{dy}{y}\right)
\tag{3.37}
$$
tends to
$$\align
&\text{Res}_{y=(q^{l_1},\ldots, q^{l_n})}
\left(\prod_{i=1}^{n}\frac{y_i^{(n-i)k}}{(1/y_i;q)_{k+1}}
\left\{(x_1y_1)^{\ld_1}\cdots (x_ny_n)^{\ld_n}+\cdots\right\}
\frac{dy}{y}\right)
\tag{3.38}\\
&=x^{\ld}~
\text{Res}_{y=(q^{l_1},\ldots, q^{l_n})}
\left(\prod_{i=1}^{n}\frac{y_i^{(n-i)k+\ld_i}}{(1/y_i;q)_{k+1}}
\frac{dy_1\cdots dy_n}{y_1\cdots y_n}\right)
+\text{lower terms}\\
&=x^{\ld}~\prod_{i=1}^{n}
\left\{\frac{(q^{-k};q)_{l_i}}{(q;q)_{k}(q;q)_{l_i}}
(q^{\ld_i+(n-i+1)k+1})^{l_i}\right\}+\text{lower terms}\,.
\endalign$$
Combining (3.32) with (3.36) and (3.38), we obtain
$$
\left(\frac{1}{2\pi\sqrt{-1}}\right)^n
\int_{\Bbb T^n} K(x;y)E_\ld(y)w(y)
\frac{dy}{y}
=C_\ld x^\ld + \text{lower terms}
\tag{3.39}
$$
in the region (3.10), 
with the leading coefficient
$$
C_\ld=\sum_{0\le l_1,\ldots, l_n\le k} \ \prod_{i=1}^{n}
\left\{\frac{(q^{-k};q)_{l_i}}{(q;q)_{k}(q;q)_{l_i}}
(q^{\ld_i+(n-i+1)k+1})^{l_i}\right\}. 
\tag{3.40}
$$
By using the $q$-binomial theorem
$$
\sum_{l \ge 0}\frac{(a;q)_l}{(q;q)_l}z^l=
\frac{(az;q)_\infty}{(z;q)_\infty} \quad (|z|<1)\,, 
\tag{3.41}
$$
the constant $C_\ld$ is determined as
$$
C_\ld=\prod_{i=1}^{n}\frac{(q^{\ld_i+(n-i)k+1};q)_k}{(q;q)_{k}}. 
\tag{3.42}
$$
This completes the proof of Proposition 3.1. 
\section{\S4: Second proof of Theorem 2.2}
In this section, we will give a proof of Theorem 2.2 
based on the principal specialization of 
nonsymmetric Macdonald polynomials, 
along the line similar to that in Sahi \cite{S}.
\par\medpagebreak
We begin with a lemma.
\Lemma{4.1}{
The function $\prod_{i=1}^n {(ux_i;q)_\infty}/{(x_i;q)_\infty}$ 
has an expansion
$$
\prod_{i=1}^n \frac{(ux_i;q)_\infty}{(x_i;q)_\infty}
=\sum_{\ld\in L^+} P_\ld(x|q,t) f_\ld(u|q,t),
\tag{4.1}
$$
in terms of Macdonald polynomials.  
The coefficients are given by
$$
f_\ld(u|q,t)=t^{n(\ld)}\prod_{s\in \ld}
\frac{1-q^{a'(s)}t^{-l'(s)}u}{1-q^{a(s)+1}t^{l(s)}}
\tag{4.2}
$$
for each partition $\ld$. 
Here, for each box $s=(i,j)$ in $\ld$, 
$a'(s)=j-1$ and $l'(s)=i-1$ stand for 
the coarm-length and the coleg-length of $s$ in $\ld$, and 
$n(\ld)=\sum_{s\in\ld}l(s)=\sum_{s\in\ld}l'(s)$. 
}
\Proof{
Let $m$ be an integer with $m\ge n$ and 
take the variables  $y=(y_1,\ldots ,y_m)$. 
Then we have 
$$
\prod_{i=1}^n
\prod_{j=1}^m
\frac{(tx_iy_j;q)_\infty}{(x_iy_j;q)_\infty}
=\sum_{\ld\in L^+} b_\ld(q,t)P_\ld(x|q,t) P_\ld(y|q,t). 
\tag{4.3}
$$
By the evaluation at $y=(1,t,\ldots ,t^{m-1})$, we get 
$$
\prod_{i=1}^n
\frac{(t^m x_i;q)_\infty}{(x_i;q)_\infty}
=\sum_{\ld\in L^+} b_\ld(q,t)P_\ld(x|q,t) P_\ld(1,t,\ldots,t^{m-1}|q,t). 
\tag{4.4}
$$
Hence we have
$$
f_\ld(t^m|q,t)=b_\ld(q,t)P_\ld(1,t,\ldots ,t^{m-1}|q,t). 
\tag{4.5}
$$
It is known by \cite{Ma2} that 
$$
P_\ld(1,t,\ldots,t^{m-1}|q,t)=
t^{n(\ld)}
\prod_{s\in\ld}
\frac{1-q^{a'(s)}t^{m-l'(s)}}{1-q^{a(s)}t^{l(s)+1}}. 
\tag{4.6}
$$
Combining this with the formula for $b_\ld(q,t)$, we obtain
$$
f_\ld(t^m|q,t)
=
t^{n(\ld)}
\prod_{s\in\ld}
\frac{1-q^{a'(s)}t^{m-l'(s)}}{1-q^{a(s)+1}t^{l(s)}}. 
\tag{4.7}
$$
Since $f_\ld(u|q,t)$ is a polynomial in $u$, and $m\ge n$ is 
arbitrary, we obtain the desired formula. 
}
We will prove Theorem 2.2 by evaluating of the function $E(x;y|q,t)$ 
at the point $y=t^\dt=(t^{n-1},t^{n-2},\ldots ,1)$.  
{}From the definition of $E(x;y|q,t)$, we easily see that 
$$
E(x;t^\dt|q,t) = \prod_{i=1}^n
\frac{(qt^nx_i;q)_\infty}{(x_i;q)_\infty}.
\tag{4.8}
$$
By Lemma 4.1, we can expand this into the sum over all $P_\ld(x|q,t)$:
$$
E(x;t^\dt|q,t) = \sum_{\ld\in L^+} 
P_{\ld}(x|q,t) f_\ld(qt^n|q,t).
\tag{4.9}
$$
For each partition $\ld\in L^+$, Macdonald's symmetric 
polynomial $P_\ld(x|q,t)$ 
can be written as a linear combination of nonsymmetric ones 
$E_\mu(x|q,t)$, summed over all compositions $\mu$ in the $W$-orbit $W.\ld$
(see \cite{Ma1}) :
$$
P_\ld(x|q,t)=\sum_{\mu\in W.\ld} a_{\ld\mu}(q,t) E_\mu(x|q,t)
\quad(a_{\ld\mu}(q,t)\in \Q(q,t))
\tag{4.10}
$$
with $a_{\ld\ld}(q,t)=1$.  
Hence we have
$$
E(x;t^\dt|q,t) = \sum_{\mu\in L} 
f_{\mu^+}(qt^n|q,t) a_{\mu^+\,\mu}(q,t) E_{\mu}(x|q,t).
\tag{4.11}
$$
On the other hand, by Theorem 2.1.(1), $E(x;y|q,t)$ has the expansion
$$
E(x;y|q,t)=\sum_{\mu\in L} a_{\mu}(q,t) E_\mu(x|q,t) E_\mu(y|q^{-1},t^{-1}).  
\tag{4.12}
$$
Hence we have
$$
E(x;t^\dt|q,t)=
\sum_{\mu\in L} a_{\mu}(q,t) E_\mu(t^\dt|q^{-1},t^{-1}) E_\mu(x|q,t). 
\tag{4.13}
$$
Comparing the two expansions (4.11) and (4.13) of $E(x;t^\dt|q,t)$, we obtain
$$
f_{\mu^+}(qt^n|q,t)a_{\mu^+\,\mu}(q,t)=
a_\mu(q,t) E_{\mu}(t^\dt|q^{-1},t^{-1}). 
\tag{4.14}
$$
In particular, if $\ld$ is a partition, then we have
$$
f_{\ld}(qt^n|q,t)=a_\ld(q,t) E_{\ld}(t^\dt|q^{-1},t^{-1}). 
\tag{4.15}
$$
\par
Evaluation of nonsymmetric Macdonald polynomials at $t^{-\dt}$ is 
already carried out by Cherednik \cite{C2}.  
If $\ld\in L^+$ is a partition, 
the value $E_\ld(t^{-\dt}|q,t)$ 
can be rewritten as follows: 
$$
E_\ld(t^{-\dt}|q,t)=t^{-(n-1)|\ld|+n(\ld)}
\prod_{s\in \ld}
\frac{1-q^{a'(s)+1}t^{n-l'(s)}}{1-q^{a(s)+1}t^{l(s)+1}}. 
\tag{4.16}
$$
{}From this formula, we have
$$
E_\ld(t^{\dt}|q^{-1},t^{-1})=t^{n(\ld)}
\prod_{s\in \ld}
\frac{1-q^{a'(s)+1}t^{n-l'(s)}}{1-q^{a(s)+1}t^{l(s)+1}}. 
\tag{4.17}
$$
\par
Substituting (4.2) and (4.17) into (4.15), we have
$$
a_\ld(q,t)=\frac{f_\ld(qt^n|q,t)}{E_\ld(t^\dt|q^{-1},t^{-1})}
=\prod_{s\in\ld}
\frac{1-q^{a(s)+1}t^{l(s)+1}}{1-q^{a(s)+1}t^{l(s)}}
\tag{4.18}
$$
for all partition $\ld\in L^+$. 
Namely we have $a_\ld(q,t)=d_\ld(q,t)/d'_\ld(q,t)$ for all $\ld\in L^+$ 
with the notation of Section 2.  
Hence by Lemma 2.5 we have
$$
a_\mu(q,t)=\frac{d_\mu(q,t)}{d'_\mu(q,t)}=
\prod_{s\in\mu}\frac{1-q^{a(s)+1}t^{l(s)+1}}{1-q^{a(s)+1}t^{l(s)}}
\tag{4.19}
$$
for all composition $\mu\in L$.
This completes the proof of Theorem 2.2.
\Remark{4.2}{
In Section 3, we determined the coefficients $a_\ld(q,t)$ by 
combining asymptotic analysis of $q$-Selberg type integrals and 
the formulas for scalar products $(E_\ld,E_\ld)$. 
Since we have derived the formulas for $a_\ld(q,t)$ 
along a different route in this section, we can also use the argument 
of Section 3 conversely to determine the scalar products 
$(E_\ld,E_\ld)$ (cf. \cite{Mi2}). 
}

\enddocument